\documentclass[12pt,leqno]{amsart}
\usepackage{amsmath,amssymb,amsfonts, amscd,   
pifont, amsthm,  amsxtra, latexsym, 
eufrak}
\usepackage[all]{xy}
\usepackage{here}
\usepackage[dvips]{graphicx}

\usepackage{psfrag}
\usepackage{amsmath}
\usepackage{amsthm}
\usepackage{amsfonts}
\usepackage{amssymb}
\usepackage[colorlinks=true,linkcolor=blue]{hyperref}
\theoremstyle{plain}
\newtheorem{thm}{Theorem}[section]
\newtheorem{theorem}[thm]{Theorem}

\newtheorem{corollary}[thm]{Corollary}

\theoremstyle{definition}
\newtheorem{definition}[thm]{Definition}
\newtheorem{prob}[thm]{Problem}
\newtheorem{example}[thm]{Example}

\newtheorem{observation}[thm]{Observation}
\newtheorem{question}[thm]{Question}

\theoremstyle{remark}
\newtheorem{remark}[thm]{Remark}

\makeatletter
\@addtoreset{equation}{section}
\def\@seccntformat#1{\csname the#1\endcsname.\quad}
\def\@seccntformat#1{\csname the#1\endcsname%
\expandafter\ifx\csname#1\endcsname\subsection.\fi\quad}
\makeatother
\newcommand{\invHom}[3]{\operatorname{Hom}_{#1}(#2,#3)}

\begin{document}
\title{Harish-Chandra's admissibility theorem and beyond}
\author{Toshiyuki Kobayashi}

\maketitle

\begin{abstract}
This article is a record of the lecture
 at the centennial conference for Harish-Chandra.  
The admissibility theorem of Harish-Chandra concerns the restrictions
 of irreducible representations to maximal compact subgroups. 
In this article, we begin with a brief explanation of two directions
 for generalizing his pioneering work
 to {\it{non-compact}} reductive subgroups:
 one emphasizes discrete decomposability
 with the finite multiplicity property, while the other focuses on finite/uniformly bounded multiplicity properties. 
We discuss how the recent representation-theoretic developments in these directions collectively offer a powerful method
 for the new spectral analysis
 of standard locally symmetric spaces, extending beyond the classical Riemannian setting.  
\end{abstract}



\newcommand{\refs}[1]{\href{#1}{\,${}^{\clubsuit\,\mbox{\scriptsize ref.}}$}}
\newcommand{\bigzero}{\Huge $0$}
\def\raisedcurvearrowright{\raisebox{.5em}{$\,\curvearrowright\,$}}


\setcounter{tocdepth}{1}
\tableofcontents
\section{Admissible Restriction {\`a} la Harish-Chandra}  
\label{sec:1}

Harish-Chandra's celebrated {\it{admissibility theorem}}
 of real reductive Lie groups $G$ 
 (Theorem \ref{thm:HCadm} below) 
 concerns the restriction of representations of $G$
 to maximal compact subgroups $K$.  
This theorem has become foundational 
 for an algebraic theory of infinite-dimensional representations
 of $G$
 using the notion of {\it{Harish-Chandra modules}}.  
It also plays a crucial role
 in proving
 that reductive groups $G$ are of type I in the Murray--von Neumann sense
 (see \cite{D77}), 
 which implies
 that irreducible decompositions of any unitary representations of $G$
 are essentially unique.

An analogous statement
 to Harish-Chandra's admissibility does not generally hold
 for the restriction to reductive symmetric pairs
 such as $(G L(n,{\mathbb{R}}), O(p,n-p))$.

In this exposition, 
 we investigate two avenues
 of research inspired by Harish-Chandra's admissibility theorem
 for the last three decades.  
We do so 
 in a more general setting, 
 focusing on the restriction 
 to {\it{non-compact}} subgroups.  
In Section \ref{sec:spec}, 
 we demonstrate 
 that these generalizations provide a new tool
 to explore the spectral theory
 of standard {\it{pseudo-Riemannian}} locally symmetric spaces.

\subsection{Reminder: Harish-Chandra's admissibility theorem}
~~~
\newline
Let $G$ be a linear real reductive Lie group, 
 and let $K$ be a maximal compact subgroup.  
We denote by $\operatorname{Irr}(K)$ the set of equivalence classes of irreducible (finite-dimensional) representations of the compact group $K$.

Harish-Chandra proved the following fundamental result \cite{HC53, HC54b}, 
 combining with Segal's theorem \cite{Se52}.

\begin{theorem}
[Harish-Chandra's admissibility theorem]
\label{thm:HCadm}
For any irreducible unitary representation $\Pi$ of $G$, 
 one has the following finite-multiplicity property:
\begin{equation}
\label{eqn:finiteK}
\text{$[\Pi|_K:\pi]<\infty$
\qquad
\text{for any }$\pi \in \operatorname{Irr}(K)$.}
\end{equation}
Here $[\Pi|_K: \pi]$ denotes the multiplicity of the representation $\pi$ 
occurring in the restriction of $\Pi$ to $K$.  
\end{theorem}

\subsection
{{}From Riemannian to reductive symmetric pairs}
~~~
\par
One may think of Harish-Chandra's admissibility
 as a theorem
 about the restriction
 of representations of $G$ to $K$, 
 symbolically written as $G \downarrow K$.  
The pair $(G,K)$ is referred to as a {\it{Riemannian symmetric pair}}
{}from a differential-geometric perspective, 
 as the homogeneous space $G/K$ carries a $G$-invariant Riemannian metric
 such that the geodesic symmetry at every point 
 defines a global isometry.  
A typical example is 
 $(G, K)=(GL(n,{\mathbb{R}}), O(n))$.  
More generally, 
 pairs such as $(G L(n, {\mathbb{R}}), O(p,q))$
 where $p+q=n$
 have a similar geometric property, 
 and are referred to as {\it{reductive symmetric pairs}}.  

\begin{definition}
\label{def:sym}
Let $G$ be  a real reductive Lie group, 
 $\sigma$ an involutive automorphism of $G$, 
 and $G'$ an open subgroup 
 of the fixed point group $G^{\sigma}$.  
The pair $(G, G')$ is called 
 a {\it{reductive symmetric pair}}.  
\end{definition}

Here are typical examples of reductive symmetric pairs.  

\begin{example}
\label{ex:sym}
{\rm{(1)}}\enspace 
The Riemannian symmetric pair $(G, K)$
 corresponding to the Cartan involution $\theta$.  
\newline
{\rm{(2)}}\enspace 
The pair $({}^{\backprime } G \times {}^{\backprime } G, \operatorname{diag}{}^{\backprime } G)$.
Let $\sigma$ be the involutive automorphism of the direct product group
 $G = {}^{\backprime } G \times {}^{\backprime } G$
 defined by $\sigma(x, y)=(y, x)$.  
Then the fixed-point subgroup $G^{\sigma} = \operatorname{diag}{}^{\backprime} G$.  
\newline
{\rm{(3)}}\enspace
The pair $(G L(n,{\mathbb{R}}), O(p,q))$
 with $p+q=n$.  
\newline
{\rm{(4)}}\enspace
The pair $(O(p,q), O(p_1,q_1) \times O(p_2, q_2))$
 with $p_1+p_2=p$, $q_1+q_2=q$.  
\end{example}

In Example \ref{ex:sym} (1), 
 the decomposition of a representation $\Pi$ of the group $G$
 with respect to the Riemannian symmetric pair $(G, K)$
 is referred to as the {\it{$K$-type formula}}.

In Example \ref{ex:sym} (2), 
 the tensor product representation 
 $\pi_1 \otimes \pi_2$
 of two representations $\pi_1$ and $\pi_2$  of the group ${}^{\backprime} G$ can be viewed as the restriction 
of the outer tensor product $\pi_1 \boxtimes \pi_2$
 of the group $G ={}^{\backprime} G \times {}^{\backprime} G$
 to the diagonal subgroup 
 $\operatorname{diag}{}^{\backprime} G$.
 
Thus, 
 the restriction of a representation of the group $G$ 
 with respect to a symmetric pair $(G, G')$
 can be considered as a generalization 
 of these two examples.

We consider two avenues
 for generalizing Harish-Chandra's admissibility theorem
 for pairs $G \supset G'$
 of reductive Lie groups, 
 particularly for reductive symmetric pairs:

$\bullet$\enspace
$G'$-admissible restriction (Section \ref{sec:admrest}), 

$\bullet$\enspace
Finite multiplicity property
 (Section \ref{sec:realsp}).

The former theme focuses on the property 
 of the irreducible decomposition of the restriction 
 $\Pi|_{G'}$ not containing continuous spectrum when the representation $\Pi$ of $G$ is unitary.  
In contrast, 
 the latter allows $\Pi$ to be non-unitary and focuses only on the finite-multiplicity property.

\vskip 1pc
\par\noindent
{\bf{Acknowledgements.}}
\newline
The article is based on the author's invited address at the 
 \textit{Centennial Conference for Harish-Chandra},
 held at Harish-Chandra Research Institute (HRI) in Prayagraj, India, from October 9 to 14, 2023.  
The author expresses his sincere gratitude
 to D.\ Prasad and A.\ Pal
 for their warm hospitality.  
This work was partially supported by the JSPS
 under the Grant-in-Aid for Scientific Research (A) (JP23H00084).

\section{Admissible Restriction to Non-Compact Subgroups}
\label{sec:admrest}
Let $G \supset G'$ be a pair of real reductive Lie groups, 
  particularly when $G'$ is a {\it{non-compact subgroup}}, 
 and $\Pi$ an irreducible unitary representation of $G$.  
In what follows, 
 we shall use the upper case letter $\Pi$
 for representations of a group $G$, 
 and the lower case letter $\pi$
 for those of the subgroup $G'$.  
Our main concern in this section 
 is how and when the restriction $\Pi|_{G'}$ behaves like 
Harish-Chandra's admissibility theorem
 (Theorem \ref{thm:HCadm}) for $\Pi|_K$.

\subsection
{Admissible restriction $G \downarrow G'$ for a non-compact subgroup}
~~~\newline
By a theorem of Mautner \cite{Mt50}, 
 the restriction $\Pi|_{G'}$ 
 of any unitary representation $\Pi$
 decomposes into a direct integral 
 of irreducible unitary representations of $G'$.  
The irreducible decomposition
 ({\it{the branching law}})
 usually contains continuous spectrum
 when $G'$ is non-compact.  

The following notion of {\it{$G'$-admissible restriction}} was introduced in \cite{xkInvent94} for a {\it{non-compact}} subgroup $G'$, partly inspired by Harish-Chandra’s admissibility theorem:

\begin{definition}
\label{def:deco}
The restriction $\Pi|_{G'}$ is said to be {\it{$G'$-admissible}}
 if it can be decomposed discretely 
into a direct sum of irreducible unitary representations
 $\pi$ of $G'$:
\[\text
{
$\Pi|_{G'} \simeq {\underset{\pi \in \widehat{G'}}\sum}^{\oplus} m_{\pi}\pi$
\quad({discrete sum})
}
\]
with the multiplicity $m_{\pi}:=[\Pi|_{G'}:\pi]$ is finite 
 for every $\pi \in \widehat{G'}$.  

Here $\widehat{G'}$ denotes the unitary dual 
 of the group $G'$, 
 that is, 
 the set of equivalence classes 
 of irreducible unitary representations of $G'$.  
\end{definition}

The key aspect of Definition \ref{def:deco} is 
 that we require not only the absence of continuous spectrum
 but also the finiteness of each multiplicity $m_{\pi}$.

When the subgroup $G'$ is compact, 
 the discrete decomposability of the restriction $\Pi|_{G'}$
 is automatically guaranteed, 
 thus the main concern is the finiteness of each multiplicity $m_{\pi}$.  
Harish-Chandra's admissibility (Theorem \ref{thm:HCadm}) corresponds to the case
 $G'=K$, 
 stating that any $\Pi \in \widehat G$
 is $K$-admissible, 
 in our terminology.

\begin{remark}
\label{rem:algdeco}
The absence of continuous spectrum is also formalized algebraically 
 in the category of $({\mathfrak{g}}, K)$-modules
 without requiring the unitarizability in \cite{xkInvent98}.  
See also Corollary \ref{cor:admrest}.  
\end{remark}

\subsection
{Restriction to compact subgroups {$K'$ $(\subset K)$}}
~~~\newline
Our interest is in analyzing the restriction $\Pi|_{G'}$
 of $\Pi \in \widehat G$
 for a pair of real reductive Lie groups $G \supset G'$.

The key idea in \cite{xkInvent94} is first to focus
 on the $K'$-structure of the $G$-module $\Pi$, 
 where $K$ and $K'$ are maximal compact subgroups
 of $G$ and $G'$, 
 respectively, 
 see below.  
\begin{alignat*}{4}
& G &&\quad \supset\quad && G'\quad
&&\text{real reductive groups}
\\[-.3em]
& \hskip 0.3pc\cup && &&\hskip 0.3pc\cup
&&
\\[-.3em]
&K &&\quad\supset\quad &&K'
&&\text{maximal compact subgroups}
\end{alignat*}

\begin{theorem}
[Criterion for $K'$-admissibility]
\label{thm:thm1}
Suppose $\Pi \in \widehat{G}$.  
Then the following two conditions (i) and (ii) 
 on the triple $(\Pi, G, K')$ are equivalent:

{\rm{(i)}}\enspace {\rm{(multiplicity)}}\enspace
$[\Pi|_{{\text{\scriptsize{$K'$}}}}:\pi]<\infty$
for any $\pi \in \operatorname{Irr}(K')$.

{\rm{(ii)}}\enspace {\rm{(geometry)}}\enspace
$\operatorname{AS}_K(\Pi) \cap C_K(K') =\{0\}$.  
\end{theorem}

The condition (ii) in Theorem \ref{thm:thm1} uses two closed cones 
 $\operatorname{AS}_K(\Pi)$ and $C_K(K')$
 in the dual space 
 of a Cartan subalgebra of ${\mathfrak{k}}$.  
$\operatorname{AS}_K(\Pi)$ is the asymptotic $K$-support
 of the representation $\Pi$.  
There are only finitely many possibilities for $\operatorname{AS}_K(\Pi)$
 among $\Pi \in \operatorname{I r r}(G)$.  
The closed cone $C_K(K')$ is the momentum set
 associated with the Hamiltonian action 
 on the cotangent bundle $T^{\ast}(K/K')$;
 see \cite[Chap.\ 6]{Ko05} for a detailed exposition.

Theorem \ref{thm:thm1} holds without assuming
 that the representation $\Pi$ is unitary.  
The implication (i) $\Rightarrow$ (ii) was proved
 in full generality 
in \cite{xkAnn98} based on an estimate of the singularity spectrum
 of hyperfunction characters
 (or the wavefront set of distribution characters).  
An alternative and algebraic proof 
 was given in \cite{K21Kostant}.  
The converse implication (ii) $\Rightarrow$ (i) was originally given in \cite{xkicm02, Ko05}
 with a sketch of the proof, 
while the full proof is provided in \cite{K21Kostant}.

\subsection{Admissible restriction $G \downarrow {G'}$}
~~~\newline
By Theorem \ref{thm:thm1}, 
 we obtain the discrete decomposability
 of the restriction of $\Pi \in \widehat G$
 with respect to $G \downarrow G'$.  

\begin{corollary}
[Criterion for admissible restriction]
\label{cor:admrest}
Suppose that a triple $(\Pi, G, G')$ 
 such that $\Pi \in \widehat G$ and $G \supset G'$ satisfies 
 one of (therefore, any of) the equivalent conditions
 in Theorem \ref{thm:thm1}.   
\newline
{\rm{(1)\enspace(\cite{xkInvent94})}}\enspace
The restriction $\Pi|_{G'}$ is {$G'$-admissible}
 (Definition \ref{def:deco}).  
\newline
{\rm{(2)\enspace(\cite{xkInvent98})}}\enspace
The underlying $({\mathfrak{g}}, K)$-module $\Pi_K$ is discretely decomposable
 as a $({\mathfrak{g}}', K')$-module 
 (Remark \ref{rem:algdeco}).  
\end{corollary}

\begin{remark}
In the case $K'=K$, 
 one has clearly $C_K(K')=\{0\}$.  
Thus Corollary \ref{cor:admrest} (1) can be viewed as a generalization 
of Harish-Chandra's admissibility (Theorem \ref{thm:HCadm}).  
\end{remark}

\subsection{Examples of admissible restrictions $\Pi|_{G'}$}
~~~\newline
For $\Pi \in \widehat G$
 and $G \supset G'$, 
 we collect typical examples of the triples $(\Pi, G, G')$
 such that the restriction $\Pi|_{G'}$ is $G'$-admissible, 
 namely, 
 it decomposes discretely with finite multiplicity, 
 as if it were Harish-Chandra's admissibility for the restriction 
 to a maximal compact subgroup $K$.  

\begin{example}
[Theta correspondence]
Let $G$ be the metaplectic group $Mp(n,{\mathbb{R}})$, 
 a double covering group
 of the symplectic group $Sp(n,{\mathbb{R}})$
 of rank $n$, 
 and $G'=G_1' \cdot G_2'$ be a compact dual pair.  
Then the restriction $\Pi|_{G'}$ 
 of the oscillator representation, 
 also referred to as the Segal--Shale--Weil representation, 
 is $G'$-admissible.  
Moreover, 
 the branching law is multiplicity-free
 {\cite{H79}}.  
\end{example}

\begin{example}
The tensor product of any two holomorphic discrete series representations
 decomposes discretely, 
 and each multiplicity is finite {\cite{R79}}.  
Moreover, 
 the multiplicities are uniformly bounded \cite{mf-korea}.  
Its generalization to reductive symmetric pairs along 
  with some explicit branching laws
 can be found 
 in \cite{mf-korea}.   
\end{example}

\subsection{Classification theory: Admissible restriction $\Pi|_{G'}$}
~~~
\par
Corollary \ref{cor:admrest} together with a necessary condition
 for the algebraic discrete decomposability
 (Remark \ref{rem:algdeco}) given in \cite{xkInvent98}
 provides a family
 of triples $(\Pi, G, G')$
 where 
$
  \Pi \in \widehat G
$ 
and 
$
G \supset G'
$
 such that the restriction $\Pi|_{G'}$ is {$G'$-admissible}, 
 namely, 
 it decomposes discretely decomposable with finite multiplicity.

Some classification results of such triples $(\Pi, G, G')$ highlighted the following cases:

$\bullet$\enspace
(\cite{KO15})\enspace
The tensor product 
$\Pi_1 \otimes \Pi_2$ for $\Pi_1$, $\Pi_2 \in \widehat G$.  

$\bullet$\enspace
(\cite{decoAq})\enspace
$(G,  G')$ is a reductive symmetric pair, 
 $\Pi$ is a discrete series representation of $G$, 
 or more generally, 
 when the underlying $({\mathfrak{g}}, K)$-module
 $\Pi_K$ is Zuckerman's derived functor module
 $A_{{\mathfrak{q}}}(\lambda)$.  

$\bullet$\enspace
(\cite{KO15})\enspace
$(G, G')$ is a reductive symmetric pair, 
 and $\Pi$ is the minimal representation of $G$.   

$\bullet$\enspace
(\cite{xdgv})\enspace
$(G, G')$ is a non-symmetric pair
 where $G'=SL(2,{\mathbb{R}})$, 
$\Pi$ is a discrete series representation of $G$.

\begin{example}
In Section \ref{sec:spec}, 
 we give yet another geometric setting
 such that the restriction $\Pi|_{G'}$ is $G'$-admissible 
 if

$\bullet$\enspace
$X$ is a $G$-space 
 on which $G'$ acts properly and spherically, 

$\bullet$\enspace
$\Pi \in \widehat G$ occurs in the space
 ${\mathcal{D}}'(X)$ of distributions.  

This geometric setting is motivated by new spectral theory
 on pseudo-Riemannian locally symmetric spaces
and the $G'$-admissibility 
is proved without relying on Theorem \ref{thm:thm1}, 
 see Theorem \ref{thm:thm10} below.  
\end{example}

For  some explicit discrete branching laws, 
 see \cite{xkInvent94, GrWa98, mf-korea, O24}; 
 for further topics
 on discretely decomposable restrictions, 
 see \cite{Ki24, KSugaku}.

\section
{Bounded/Finite Multiplicity Pairs for Restriction}
\label{sec:realsp}
Let $G \supset G'$ be a pair
 of real reductive Lie groups.  
In the previous section,
 we focused on $G'$-admissible restrictions
 of irreducible unitary representations $\Pi$, 
 as a way to generalize Harish-Chandra's admissibility 
 to non-compact subgroups
 by requiring discrete decomposability
 and finiteness of multiplicities.  
The absence of continuous spectrum facilitates an algebraic approach
 to the restriction $\Pi|_{G'}$
 of the unitary representation $\Pi$, 
 even when $G'$ is non-compact.

In this section, 
 we focus solely 
 on the multiplicity 
 by dropping the requirement
 of discrete decomposability.

We will see
 that the finite-multiplicity property is not obvious 
 even for reductive symmetric pairs, 
 such as 
\[
\text{$
   (G L(n,{\mathbb{R}}), O(p,q))
$
 or 
$
   (G L(n,{\mathbb{R}}), G L(p,{\mathbb{R}}) \times G L(q,{\mathbb{R}}))
$
 when $p+q=n$.}
\]  
We also explore a stronger condition 
 referred to as 
 {\it{bounded multiplicity property}},
 which does not generally hold even 
 for restrictions related to Riemannian symmetric pairs
 $(G, K)$
 but it still holds for some reductive symmetric pairs
 (Theorem \ref{thm:thm3}).  
These perspectives lead us 
 to yet another avenue of the restriction problem.

\subsection{Reminder: smooth representations of moderate growth}
\label{subsec:admrep}
~~~\newline
It is observed 
 that irreducible continuous representations $\Pi$
 of real reductive groups can exhibit wild behavior 
 if they are not unitary.  
Even when $G={\mathbb{R}}$, 
  there exists an infinite-dimensional irreducible
 representation
 of the abelian group $G$
 on a Banach space, 
 as a consequence of a counterexample
 to the invariant subspace problem by Enflo, {\it{cf.}} \cite{E87}.  
Harish-Chandra's admissibility 
 provides a guiding principle 
 for identifying \lq\lq{reasonable}\rq\rq\
 classes of continuous representations
 of reductive Lie groups
 avoiding such counterexamples, 
 defined as below.  

\begin{definition}
[admissible representation]
A continuous representation $\Pi$ of $G$ is said to be {\it{admissible}} 
({\it{$K$-admissible}} for later terminology)
 if 
\[\text
{$[\Pi|_K:\pi]<\infty$
\qquad
for every $\pi \in \operatorname{Irr}(K)$.}
\]
\end{definition}

To be rigorous about \lq{multiplicities}\rq\
 for infinite-dimensional representations, 
 we need to specify the topology 
 of the representation spaces.  
For this purpose, 
 let $G$ be a real reductive Lie group, 
 and ${\mathcal{M}}(G)$ denote the category 
 of smooth admissible representations
 of finite length with moderate growth, 
 which are defined on Fr{\'e}chet topological vector spaces \cite[Chap.\ 11]{WaI}.  
Casselman--Wallach's theory shows
 that there is a natural category equivalence 
 between ${\mathcal{M}}(G)$
and the category of $({\mathfrak{g}}, K)$-modules 
 of finite length.  

Let $\operatorname{I r r}(G)$ denote 
 the set of irreducible objects in ${\mathcal{M}}(G)$.

If $\Pi$ is an admissible continuous representation
 of finite length on a Banach space, 
 then the representation $\Pi^{\infty}$
 acting on the Fr{\'e}chet space
 of $C^{\infty}$-vectors
 belongs to ${\mathcal{M}}(G)$.  
In particular, 
 this yields a natural injection:
\[
  \underset{\text{\scriptsize{{unitary dual}}}}{\widehat G} \hookrightarrow \operatorname{Irr}(G), 
\quad
 \Pi \mapsto \Pi^{\infty}.  
\]

\vskip 1pc
\subsection
{Multiplicity of the restriction $\Pi|_{G'}$}
~~~
\par

We use the category ${\mathcal{M}}(G)$
 to define the \lq\lq{multiplicity}\rq\rq\
 in the restriction $\Pi|_{G'}$.  

\begin{definition}
[Symmetry Breaking Operator]
\label{def:SBO}
Let $\Pi \in {\mathcal{M}}(G)$
 and $\pi\in {\mathcal{M}}(G')$.  
A continuous $G'$-homomorphism from $\Pi|_{G'}$ to $\pi$ is referred
 to as a {\it{symmetry breaking operator}}.  
Let $\operatorname{Hom}_{G'}(\Pi|_{G'}, \pi)$ denote 
 the vector space of symmetry breaking operators.  
The {\it{multiplicity}} of $\pi$ in the 
restriction $\Pi|_{G'}$ is defined
 by its dimension, 
 that is, 
\begin{equation}
\label{eqn:mult}
[\Pi|_{G'}:\pi]
:=
\dim
{\operatorname{Hom}_{G'}
(\Pi|_{G'}, \pi)}
\in {\mathbb{N}} \cup \{\infty\}.  
\end{equation}
\end{definition}

The definition of the multiplicity \eqref{eqn:mult}
 in the category ${\mathcal{M}}(G)$
 coincides with the multiplicity 
in the category of unitary representations
 if the restriction is $G'$-admissible
 (Definition \ref{def:deco}).

\subsection
{Comparison: $GL(n,{\mathbb{R}}) \downarrow O(n)$\,\, 
{\it{vs}}\,\, $GL(n,{\mathbb{R}}) \downarrow O(p,n-p)$}
~~~\newline
Harish-Chandra's admissibility theorem 
 (Theorem \ref{thm:HCadm}) concerns
 the restriction
 with respect to a {Riemannian symmetric pair}
\[
G \supset K, 
\qquad
\text{{e.g., {$GL(n,{\mathbb{R}}) \supset O(n)$}}}
\]
and asserts that 
\[\text
{$
   [\Pi|_K: \pi]<\infty
\quad
\text{for any \,$\Pi \in \operatorname{Irr}(G)$\, and for any \,$\pi \in \operatorname{Irr}(K)$.}
$}
\]

In contrast, 
 for a reductive symmetric pair
\[
\text{
$G \supset G'$, 
{\it{e.g.,}} 
$
GL(n,{\mathbb{R}}) \supset O(p,n-p)
$}
\]
it may occur that 
\[
[\Pi|_{G'}: \pi]=\infty
\quad
\text{for some $\Pi \in \operatorname{Irr}(G)$
 and $\pi \in \operatorname{Irr}(G')$}.  
\]

The classification of reductive symmetric pairs $(G, G')$
 having the finite multiplicity property
\[
  [\Pi|_K: \pi]< \infty\quad
  \text{for any $\Pi \in \operatorname{I r r}(G)$
 and for any $\pi \in \operatorname{I r r}(G')$}.  
\]
has been established in \cite{xKMt}
 using the criterion
 in Theorem \ref{thm:thm4} below.  
We will explain the background of the theory.

\subsection{Spherical vs real spherical}
~~~\newline
Let $G_{\mathbb{C}}$ be a complex reductive Lie group, 
 and $X_{\mathbb{C}}$ a connected complex manifold 
 on which $G_{\mathbb{C}}$ acts holomorphically.

\begin{definition}
\label{def:spherical}
The $G_{\mathbb{C}}$-space $X_{\mathbb{C}}$ is {\it{spherical}} if
a Borel subgroup $B$ of $G_{\mathbb{C}}$ has an open orbit in $X_{\mathbb{C}}$.
\end{definition}

\begin{example}
[{\cite{Wo69}}]
Complex reductive symmetric spaces are spherical.  
\end{example}

In search of a {broader framework for global analysis}
 on homogeneous spaces than those known
 in the late 1980s 
 (e.g. group manifolds, semisimple symmetric spaces), 
 the author advocated
 introducing the following concept from the viewpoint
 of the finite multiplicity property.

\begin{definition} 
[{\cite{K95}}]
\label{def:Rspherical}
Let $X$ be a connected $C^{\infty}$ manifold
 on which a real reductive Lie group $G$ acts as diffeomorphisms.  
We say $X$ is {\it{real spherical}}
 if a minimal parabolic $P$ of $G$ has an open orbit in $X$.  
\end{definition}

In what follows, 
 we assume that $G$ is a real reductive group
 with complexification $G_{\mathbb{C}}$, 
 and that $H$ is an algebraic reductive subgroup of $G$
 with complexification $H_{\mathbb{C}}$.  
We write $X=G/H$ and $X_{\mathbb{C}}=G_{\mathbb{C}}/H_{\mathbb{C}}$.

\begin{remark}
It is important 
 to emphasize that the following two notions differ:

$\bullet$\enspace real forms of spherical spaces, 

$\bullet$\enspace real spherical spaces.  

The former is stronger than the latter;
 namely, 
 if $G_{\mathbb{C}}/H_{\mathbb{C}}$ is spherical
 then its real form $G/H$ is always real spherical
 (\cite[Lem.\ 3.5]{xktoshima}).  
The converse is not necessarily true.  
For instance, 
 any homogeneous $G/H$ is real spherical
 if $G$ is compact, 
 but its complexification $G_{\mathbb{C}}/H_{\mathbb{C}}$
 is not necessarily spherical.   
\end{remark}

\begin{remark}
{\rm{(1)}}
The definition of real sphericity of $X=G/H$
 in Definition \ref{def:Rspherical}
 is equivalent to the condition 
 that $\#(P\backslash G/H)<\infty$ (Kimelfeld, Matsuki, Bien).  
This extends a theorem by Brion and Vinberg, 
which asserts that a $G_{\mathbb{C}}$-space $X_{\mathbb{C}}$
 is spherical if and only if 
 $\# (B \backslash X_{\mathbb{C}}) < \infty$.   
\newline
{\rm{(2)}}\enspace
For compact $H$, 
Akhiezer--Vinberg \cite{AV99} proved
 that $G/H$ is a weakly symmetric space in the sense of Selberg
 if and only if $G_{\mathbb{C}}/H_{\mathbb{C}}$ is spherical.  
\newline
{\rm{(3)}}\enspace
For compact $G$, 
 Tanaka \cite{TanakaJDG22} proved that $X_{\mathbb{C}}$ is $G$-strongly visible
 in the sense of \cite{Ko05}
 if and only if $X_{\mathbb{C}}$ is $G_{\mathbb{C}}$-spherical.  
\end{remark}

\subsection{Analytic view of real spherical spaces}
~~~\newline
A fundamental requirement
 in non-commutative harmonic analysis
 for a $G$-space $X$ is 
 that the space of functions on $X$ should be well-controlled 
 by the group $G$.  
To be rigorous, 
 we formalize the degree of control 
 of the group $G$
 in terms of the multiplicity
 $\dim \invHom G \Pi {C^{\infty}(X)}$, 
 that is, 
 the number of times each irreducible representation
 $\Pi \in \operatorname{I r r}(G)$
 occurs in $C^{\infty}(X)$.

The following theorem provides a geometric criterion for the finiteness
 of the multiplicity:

\begin{theorem}
[Finite Multiplicity Space]
\label{thm:rs}
Let $G$ be a reductive Lie group, 
 $H$ a reductive algebraic subgroup of $G$, 
 and $X=G/H$.  
Then the following two conditions on the pair $(G,H)$ are equivalent.  
\newline
{\rm{(i)}}\enspace
{\rm{(Global analysis)}}\enspace
$\dim \invHom G \Pi{C^{\infty}(X)}< \infty$
for every $\Pi \in \operatorname{I r r}(G)$.  
\newline
{\rm{(ii)}}\enspace
{\rm{(Geometry)}}\enspace
$X$ is $G$-real spherical.  
\end{theorem}

In the original proof (\cite[Thm.\ A]{xktoshima}), 
 we consider a more general setting 
 for the space of distribution sections 
 of equivariant vector bundles
 over $X=G/H$, 
 where $H$ is not necessarily reductive, 
 and not only gives a qualitative result
 (the equivalence in Theorem \ref{thm:rs})
 but also provides quantitative results.  
Specifically, 
 we give an upper estimate of the multiplicity
 by using hyperfunction boundary value maps
 for partial differential equations.  
For a lower estimate, 
 we generalize an idea of the classical Poisson transform, 
 see \cite[Sec.\ 6.1]{K14} for more details.  
These estimates from the above and below
 establish a necessary and sufficient condition
 for the uniform boundedness 
 of the multiplicity, 
 which gives a stronger degree of control
 of the group $G$
 over the function space of $X$, 
 as discussed in the next section.


\subsection{Analytic view of spherical spaces}
~~~
\newline
A discovery in \cite{xrims40, xktoshima}
 reveals the fact 
 that the uniform boundedness property of the multiplicity
 in $C^{\infty}(X)$ is determined
 solely by the complexification $X_{\mathbb{C}}=G_{\mathbb{C}}/H_{\mathbb{C}}$.
This is in sharp contrast
 to the finiteness criterion 
 established in Theorem \ref{thm:rs}.

The results are summarized 
 in the following theorem, 
 which shows a harmony 
 of analysis, geometry, and algebra:

\vskip 1pc
\[
\begin{array}{ccc}
{\text{Geometry}} 
& & 
{\text{Analysis}}
\\[0.5em]
G_{\mathbb{C}} \raisedcurvearrowright X_{\mathbb{C}} 
&
\rightsquigarrow
&
G \raisedcurvearrowright C^{\infty}(X)
\\[0.3em]
\hphantom{MMMM}
\rotatebox[origin=c]{315}{$\rightsquigarrow$}
&
&\rotatebox[origin=c]{45}{$\rightsquigarrow$}
\hphantom{MMMM}
\\[0.3em]
&{\mathbb{D}}_G(X)
&
\\[0.3em]
&{\text{Algebra}} 
&
\end{array}
\]

\begin{theorem}
[Criterion for Uniformly Bounded Multiplicity]
\label{thm:thm2}
Let $G \supset H$ be a pair of real reductive Lie groups, 
 and $X=G/H$.  
Then the following conditions (i), (ii), (iii), and (iii)$'$
 on the pair $(G, H)$
 are equivalent:
\newline
{\rm{(i) (Global analysis)}}\enspace
There exists a constant $C>0$ such that
\begin{equation}
\label{eqn:bddX}
 \dim \operatorname{Hom}_G (\pi, C^\infty (X)) \le C 
 \quad\text{for all $\pi \in \operatorname{Irr}(G)$}.
\end{equation}
\newline
{\rm{(ii) (Complex geometry)}}\enspace 
$X_{\mathbb{C}}$ is $G_{\mathbb{C}}$-spherical.
\newline
{\rm{(iii) (Ring structure)}}\enspace
The algebra ${\mathbb{D}}_{G_{\mathbb{C}}}(X_{\mathbb{C}})$ is a commutative ring. 
\newline
{\rm{(iii) $'$(Ring structure)}}
The algebra ${\mathbb{D}}_{G_{\mathbb{C}}}(X_{\mathbb{C}})$ is
 a polynomial ring.  
\end{theorem}

The proof of the equivalence (iii) $\Leftrightarrow$ (iii)$'$ is given
 in Knop \cite{Kn94}, 
 and the equivalence (ii) $\Leftrightarrow$ (iii) was established earlier, 
 see Vinberg \cite{Vi01}
 and references therein.

In contrast to conditions (ii), (iii), and (iii)$'$, 
 depending solely on the complexifications
 $(G_{\mathbb{C}}, H_{\mathbb{C}})$, 
 the objects in (i) such as $\operatorname{I r r}(G)$
 and $C^{\infty}(G)$ are strongly dependent
 on the choice of the real forms $G$ and $X$ of $G_{\mathbb{C}}$
 and $X_{\mathbb{C}}$.  
The equivalence (i) $\Leftrightarrow$ (ii) 
 was established by the author
 in collaboration with T.\ Oshima \cite{xktoshima}.

\subsection
{Restriction $G\downarrow G'$ with {finite multiplicity property}}
~~~\newline
We apply Theorem \ref{thm:rs} 
 to the homogeneous space 
 $(G \times G')/\operatorname{diag}G' \simeq G$
 to study the restriction of representations 
 with respect to $G \downarrow G'$.  

\begin{theorem}
[{Finite Multiplicity Pairs for Restriction, {\cite{K14}}}]
\label{thm:thm4}
The following two conditions (i) and (ii) for a pair 
 of real reductive groups $G \supset G'$
 are equivalent:
\newline
{\rm{(i)}}\enspace{\rm{(Representation theory)}}\enspace
$[\Pi|_{G'}:\pi]<\infty$
\quad
for every $\Pi \in \operatorname{Irr}(G)$, 
and for every $\pi \in \operatorname{Irr}(G')$.  
\newline
{\rm{(ii)}}\enspace{\rm{(Geometry)}}\enspace
$(G \times G') / \operatorname{diag}G'$
 is real spherical.  
\end{theorem}

\begin{example}
[Harish-Chandra's admissibility]
The geometric condition (ii) in Theorem \ref{thm:thm4} holds when $G'=K$
 by the Gauss--Iwasawa decomposition
 $G=K A N$.  
The representation-theoretic condition (i) corresponds to 
 Harish-Chandra's admissibility when $G'=K$.  
\end{example}

A complete classification of reductive
 symmetric pairs
 $(G,G')$, 
 where $G'$ is {\it{non-compact}}, 
 satisfying the geometric condition (ii)
 was established by the author
 in collaboration with Matsuki \cite{xKMt}.

\subsection
{Restriction $G\downarrow G'$
 with {uniformly bounded multiplicity property}}
~~~\newline
As in Theorem~\ref{thm:thm4}, applying Theorem~\ref{thm:thm2} to the homogeneous space $(G \times G')/\operatorname{diag}G'$ yields a necessary and sufficient condition for the uniform boundedness property of the restriction $G \downarrow G'$:

\begin{theorem}
[Uniformly Bounded Multiplicity Pairs for Restriction]
\label{thm:thm3}
Let {$G \supset G'$} be a pair 
 of real reductive groups.  
Then the following four conditions (i), (ii), (iii), and (iii)$'$
 are equivalent:

{\rm{(i)}}\enspace\hphantom{ii}
{\rm{(Representation Theory)}}\enspace
$\underset{\Pi \in \operatorname{Irr}(G)}\sup\,\,
\underset{\pi \in \operatorname{Irr}(G')}\sup\,\,
[\Pi|_{G'}:\pi]<\infty$.

{\rm{(ii)}}\enspace\hphantom{i}
{\rm{({Complex Geometry})}}\enspace
$(G_{\mathbb{C}} \times G_{\mathbb{C}}') / \operatorname{diag}G_{\mathbb{C}}'$
 is spherical.  

{\rm{(iii)}}\enspace\hskip 0.2pc
{\rm{({Ring})}}\enspace
The algebra $U({\mathfrak{g}}_{\mathbb{C}})^{G'}$ is a commutative ring.  

{\rm{(iii)$'$}}\enspace{\rm{({Ring})}}\enspace
The algebra $U({\mathfrak{g}}_{\mathbb{C}})^{G'}$ is a polynomial ring.  
\end{theorem}

See \cite{K14} 
 for some other equivalent conditions, 
 as well as for the proof of the equivalence (i) $\Leftrightarrow$ (ii).  
See also \cite{K95, xktoshima}.  
There are a few pairs $(G, G')$
 satisfying (i)--(iii)$'$
 but in which the supremum in (i) is greater than one.  
However, 
 for most of the important cases, 
 a sharper estimate
 for (ii) $\Rightarrow$ (i) holds, 
 that is, the supremum in (i) equals one 
 (Sun--Zhu \cite{xSunZhu12}).

In the trivial case where $G=G'$, 
 the finiteness condition (i) is evident.  
The sphericity condition (ii) is guaranteed 
 by the Bruhat decomposition, 
 while the ring structure (iii)$'$ is established
 via the Harish-Chandra isomorphism \cite{HC8}.

The classification of such complex pairs
 $({\mathfrak{g}}_{\mathbb{C}},{\mathfrak{g}}_{\mathbb{C}}')$ was provided
 by Kr{\"a}mer \cite{Kr76} and Kostant in 1970s
 if ${\mathfrak{g}}_{\mathbb{C}}$ is simple, 
 specifically 
 $({\mathfrak{g}}_{\mathbb{C}},{\mathfrak{g}}_{\mathbb{C}}')$ being 
 $({\mathfrak{sl}}(n,{\mathbb{C}}),{\mathfrak{gl}}(n-1,{\mathbb{C}}))$, 
$({\mathfrak{so}}(n,{\mathbb{C}}),{\mathfrak{so}}(n-1,{\mathbb{C}}))$, 
 up to abelian factors
 and possibly considering outer automorphisms.

\subsection
{Example: $O(p,q) \downarrow O(p_1,q_1) \times O(p_2,q_2)$}
~~~\newline
We examine the finiteness criterion (Theorem \ref{thm:thm4})
 and the uniform boundedness criterion 
 (Theorem \ref{thm:thm3}) in the context
 of the symmetric pair $(G,G')=(O(p,q), O(p_1,q_1) \times O(p_2,q_2))$
 where $p_1+p_2=p$
 and $q_1 + q_2=q$, 
 and $p+q \ge 5$.

Our criteria tell us the following equivalences:

$\bullet$\enspace
$\underset{\Pi \in \operatorname{Irr}(G)}{\operatorname{sup}}\,\,
 \underset{\pi \in \operatorname{Irr}(G')}{\operatorname{sup}}
 [\Pi|_{G'}:\pi]<\infty
$
$\,\,\iff\,\,$
$
\text{$p_1+q_1=1$ or $p_2+q_2=1$;}
$

$\bullet$\enspace
$[\Pi|_{G'}:\pi]<\infty$
\quad
for every $\Pi \in \operatorname{Irr}(G)$ 
 and for every $\pi \in \operatorname{Irr}(G')$

\par\indent
$\hphantom{MMM}\,\,{\iff}\,\,$
$p_1+q_1=1$, $p_2+q_2=1$, 
 $p=1$, $q=1$, or $G'$ compact.

The proof of Theorem \ref{thm:thm4} shows that 
 for general values of $p_1$, $p_2$, $q_1$, $q_2$, 
 it can occur that the multiplicity
\[
  [\Pi|_{G'}:\pi]=\infty
\]
for some principal series representations $\Pi \in \operatorname{Irr}(G)$
 and $\pi \in \operatorname{Irr}(G')$, 
 which contrasts with the case
 where $G'$ is compact 
({Harish-Chandra's admissibility}).

\subsection
{Question: Bounded multiplicity $\Pi|_{G'}$ for \lq\lq{small}\rq\rq\ $\Pi$}
~~~\newline
By inspecting the above examples, 
 one sees
 that refining the question should broaden 
 the concept of 
 \lq\lq{good classes}\rq\rq\
 for branching problems.  
Thus,
 we consider a triple $(\Pi, G, G')$, 
 where $\Pi \in \operatorname{I r r}(G)$
 and $G \supset G'$, 
 instead of just a pair $(G, G')$, 
 as in the formulation of the admissible restriction 
 (Corollary \ref{cor:admrest}).  
We now pose the following question.  
\begin{question}
\label{q:Piuni}
Given a reductive symmetric pair {$G \supset G'$}, 
 does there exist {at least one}
 infinite-dimensional $\Pi \in \operatorname{Irr}(G)$
with the following bounded multiplicity property?

\[\text{
$\underset{\pi \in \operatorname{Irr}(G')}{\operatorname{sup}}
[\Pi|_{G'}:\pi]<\infty$.
}\]  
\end{question}

\subsection{An affirmative answer to Question \ref{q:Piuni}}
~~~\newline
In \cite{mf-korea},
 we provided an affirmative answer to Question \ref{q:Piuni}
in the context of the reductive symmetric pair $(G, G')$
 given by Hermitian Lie groups such as $(U(p,q), U(p_1,q_1) \times U(p_2,q_2))$
 where $p_1+p_2=p$ and $q_1 +q_2 =q$
 using the theory of {\it{visible actions}}
 on complex manifolds.  
For the general case, 
 we have proved the following result.  

\begin{theorem}
[{\cite{K22}}]
\label{thm:thm5}
Let $G$ be a simply connected, non-compact, real semisimple Lie group.  
There exist a constant $C\equiv C(G) >0$ and an
 infinite-dimensional irreducible representation $\Pi$ of $G$
 such that 
\[\text{
$\underset{\pi \in \operatorname{Irr}(G')}{\operatorname{sup}}
[\Pi|_{G'}:\pi]\le C$
}\]
for {\it{all}} symmetric pairs $(G, G')$.  
\end{theorem}

\begin{example}
[tensor product]
There exist a constant $C>0$ and
 infinite-dimensional 
 irreducible representations $\Pi_1$, $\Pi_2$ of $G$
 such that 
\[\text{
$[\Pi_1 \otimes \Pi_2:\Pi] \le C$}
\quad
\text{for every $\Pi \in \operatorname{Irr}(G)$}.  
\]
\end{example}

It is important to note that Theorem~\ref{thm:thm5} establishes the uniformly bounded multiplicity property.
Thus, even in the special case of $G' = K$ in Theorem~\ref{thm:thm5}, the result does not follow from Harish-Chandra’s admissibility theorem, which guarantees only individual finiteness but not the uniformly bounded multiplicity property.

Another example of $\Pi \in \widehat G$ that gives an affirmative answer to Question~\ref{q:Piuni} is the minimal representation:
\begin{theorem}
[{\cite{tkVarnaMin}}]
\label{thm:thm6}
Let $G$ be a real reductive Lie group.  
If the associated variety
 of $\Pi \in \operatorname{Irr}(G)$
 is the minimal nilpotent orbit 
 in ${\mathfrak{g}}_{\mathbb{C}}^{\ast}$, 
then there exists a constant $C \equiv C(\Pi) >0$
 such that
\[
\text
{
$
\underset{\pi \in \operatorname{Irr}(G')}{\sup}[\Pi|_{G'}:\pi]\le C
$
}
\]
for {all} reductive symmetric pairs $(G,G')$.  
\end{theorem}

\subsection{Methods of proof}
~~~
\newline
The original approach in \cite{xktoshima}
 to prove the finiteness/uniform boundedness
 for multiplicities 
under certain geometric conditions
 ({\it{e.g.}}, (ii) $\Rightarrow$ (i) in Theorem \ref{thm:rs}, 
(ii) $\Rightarrow$ (i) in Theorem \ref{thm:thm2}, 
 and (ii) $\Rightarrow$ (i) in Theorem \ref{thm:thm4}, 
etc.)
was to use partial differential equations
 with regular singularities
 at the boundary of a specific compactification.  
Further approaches used 
 by Tauchi, Kitagawa, 
  Kobayashi, Aisenbud--Gourevich, Tanaka and other researchers
 include holonomic ${\mathcal{D}}$-modules \cite{Tu}, 
 visible actions on complex manifolds \cite{xrims40, Tn24},
 and coisotropic actions on symplectic manifolds
 \cite{Ki}, etc.

Some of these methods have broader applications, 
 though upper estimates
 of the multiplicities
 are not necessarily 
 as sharp as those in \cite{xktoshima}.

The proof of the converse statement employs
 integral transforms \cite{K14}, 
 which provide lower bounds
 for the multiplicities in Theorems \ref{thm:thm2}
 and \ref{thm:thm3}.

\section{Bounded Multiplicity Triple for Restriction}
\label{sec:small}

In search of a natural framework for a detailed and potentially fruitful analysis of branching problems (for example, Stage C in the ABC program \cite{xKVogan2015}, which studies the restriction of representations), we focus on a specific family of “small” representations $\Pi \in \operatorname{Irr}(G)$ for which the restriction has the {\it{uniformly bounded}} property:
By {\it{uniformity}}, 
 we consider not only representations {$\pi \in \operatorname{Irr}(G')$}
 but also the family of representations $\Pi \in \operatorname{Irr}(G)$
 as described below.

\begin{question}
\label{q:supsup}
Find a triple $(\Omega, G, G')$
 where {$\Omega \subset \operatorname{Irr}(G)$}
 and $G \supset G'$
 which satisfies 
 the following uniform boundedness property:
\begin{equation}
\label{eqn:supsup}
\text
{$\underset{{\text{\scriptsize{$\Pi \in \Omega$}}}}{\operatorname{sup}}\,\,
\underset{{\text{\scriptsize{$\pi \in \operatorname{Irr}(G')$}}}}{\operatorname{sup}}
[\Pi|_{G'}:\pi] < \infty$.  
}
\end{equation}
\end{question}

This question further explores Question \ref{q:Piuni}.  
We have seen several affirmative results 
 in the previous sections, 
 such as
\begin{alignat*}{5}
&\,\,\Omega=\operatorname{Irr}(G)
&&
&&
&&
&&
\text{(Theorem \ref{thm:thm3})}
\\[.8em]
&\,\,\Omega= \{\text{minimal representations}\}
\qquad
&&
&&
&&
&&\text{(Theorem \ref{thm:thm6})}
\end{alignat*}

In the next section, 
 we will discuss Question \ref{q:supsup}, 
 focusing on $\Omega:=\{\text{$H$-distinguished representations}\}$.

\subsection
{Restriction of $H$-distinguished representations $H \nearrow G$}
~~~\newline
We set up some notation.  
Let $H$ be a closed subgroup of a Lie group $G$.  

\begin{definition}
\label{def:distinguished}
We say $\Pi\in \operatorname{Irr}(G)$
 is an {\it{$H$-distinguished representation}}
 of $G$, 
 if {$(\Pi^{-\infty})^H \ne \{0\}$}, 
 or equivalently if 
\begin{equation*}
\operatorname{Hom}_G(\Pi,\, C^{\infty}(G/H)) \ne \{0\} 
\end{equation*}
by the Frobenius reciprocity.  
Let $\operatorname{Irr}(G)_H$ denote 
 the subset of $\operatorname{Irr}(G)$
 consisting
 of $H$-distinguished irreducible admissible representations.
\end{definition}

\subsection
{Bounded multiplicity triple for $H \nearrow G \searrow G'$}

\begin{definition}
\label{def:3.3}
A triple $H \subset G \supset G'$ of real reductive Lie groups
 is said to be a {\it{bounded multiplicity triple}}
 if \eqref{eqn:supsup} holds
 for $\Omega= \operatorname{I r r}(G)_H$.  
\end{definition}

This means
 that for every representation $\Pi \in \operatorname{Irr}(G)$
 that can be realized in $C^{\infty}(G/H)$, 
 the multiplicity of the restriction $\Pi|_{G'}$
 is uniformly bounded:
\[
   \underset{{\text{\scriptsize{$\Pi \in \operatorname{I r r}(G)_H$}}}}{\operatorname{sup}}\,\,
\underset{{\text{\scriptsize{$\pi \in \operatorname{Irr}(G')$}}}}{\operatorname{sup}}
[\Pi|_{G'}:\pi] < \infty.  
\]

We shall discuss an aspect of the classification theory
 of bounded multiplicity triples 
 established in \cite{K22} below.

We begin with some na{\"i}ve considerations for Definition \ref{def:3.3}
 including:

(1)\enspace
If the subgroup $H$ is \lq\lq{large}\rq\rq\ in $G$, 
 then we may expect that the representation $\Pi \in \operatorname{Irr}(G)_H$
 will be \lq\lq{small}\rq\rq.

(2)\enspace
If the subgroup $G'$ is \lq\lq{large}\rq\rq\ in $G$
 and if the representation $\Pi$ is \lq\lq{small}\rq\rq, 
 then we may expect that the subgroup $G'$ will have 
 a strong degree of the control, 
 in the sense that the restriction $\Pi|_{G'}$
 will have a bounded multiplicity property.  
\vskip 0.5pc
The \lq\lq{largeness}\rq\rq\
 of the subgroups
 depends on the properties
 we are examining.  
The conditions for $H$ are formulated
 in Theorem \ref{thm:GHquotient}, 
 while those for $G'$ and $H$ are presented
 in Theorem \ref{thm:thm7} (ii) below.

\subsection
{Relative Borel subalgebra ${\mathfrak{b}}_{G/H}$
 and relative parabolic subgroup $P_{G/H}$}
~~~\newline
In non-commutative harmonic analysis of real reductive Lie groups $G$, 
 the notion of minimal parabolic subgroups of $G$ and Borel subalgebras
 in the complexified Lie algebras ${\mathfrak{g}}_{\mathbb{C}}$ 
plays a fundamental role.  
We consider the generalization
 of this notion  to reductive symmetric spaces
 $G/H$ 
 associated to an involutive automorphism $\sigma$ of $G$.

We take a maximal compact subgroup $G_U \subset G_{\mathbb{C}}$
 such that $G_U \cap G$
 and $G_U \cap H$ are also maximal compact subgroups 
 of $G$ and $H$, 
 respectively.  
We fix an $\operatorname{Ad}(G)$-invariant non-degenerate symmetric bilinear form
 on ${\mathfrak{g}}$
 which is also non-degenerate on the subalgebra ${\mathfrak{h}}$.  
We write ${\mathfrak{g}}={\mathfrak{h}}+{\mathfrak{h}}^{\perp}$
 for the direct sum decomposition
 and ${\mathfrak{g}}_{\mathbb{C}}={\mathfrak{h}}_{\mathbb{C}}+{\mathfrak{h}}_{\mathbb{C}}^{\perp}$
 for its complexification.  
Recall that to a given hyperbolic element $Y$ in ${\mathfrak{g}}$, 
 one associates a parabolic subalgebra of ${\mathfrak{g}}$, 
 defined as the sum of eigenspaces
 of $\operatorname{a d}(Y)$
 with non-negative eigenvalues.  
\begin{definition}
\label{def:PGH}
(Relative Borel subalgebra ${\mathfrak{b}}_{G/H}$  
 and parabolic
 subalgebra ${\mathfrak{p}}_{G/H}$ {\cite{K22}}).  
Let $(G,H)$ be a reductive symmetric pair.  
\newline
(1)\enspace
 {{\it{A Borel subalgebra} {${\mathfrak{b}}_{G/H}$}}} for $G/H$
is a parabolic subalgebra
 of ${\mathfrak{g}}_{\mathbb{C}}$.  
It is defined 
 by a generic element of 
$
   {\mathfrak{h}}_{\mathbb{C}}^{\perp} \cap \sqrt{-1}{\mathfrak{g}}_U
$
 or its conjugate
 by an inner automorphism of $G_{\mathbb{C}}$.  
\newline
(2)\enspace
{{\it{A minimal parabolic subalgebra}} {${\mathfrak{p}}_{G/H}$} for $G/H$}
 is a real
 parabolic subalgebra of ${\mathfrak{g}}$.  
It is defined by a generic 
element
 of ${\mathfrak{h}}^{\perp} \cap \sqrt{-1} {\mathfrak{g}}_U$
 or its conjugate
 by an inner automorphism of $G$. 
\end{definition}

According to the definition, 
${\mathfrak{b}}_{G/H}$ is determined 
 solely from the complexified symmetric pair $({\mathfrak{g}}_{\mathbb{C}}, {\mathfrak{h}}_{\mathbb{C}})$.

\begin{remark}
In contrast to the usual definition of a Borel subalgebra, 
 the relative Borel subalgebra {${\mathfrak{b}}_{G/H}$} is not necessarily solvable.  
\end{remark}

\subsection
{Bounded multiplicity theorem for $H \nearrow G \searrow G'$}
~~~\newline
In Theorems \ref{thm:thm7} and \ref{thm:thm8} below, 
 we set the following conditions:
\begin{alignat*}{2}
 &(G, H) \quad 
&&\text{: a reductive symmetric pair, }
\\
&G'
&&\text{: a reductive subgroup of $G$.}
\end{alignat*}

We do not need to assume
 that $(G,G')$ is a symmetric pair.  

\begin{theorem}
[Bounded Multiplicity Criterion, {\cite[Thm.\ 1.2]{K22}}]
\label{thm:thm7}
~~~\newline
The following two conditions 
 {\rm{(i)}} and {\rm{(ii)}} on a triple $H \subset G \supset G'$
 are equivalent.  
\\
{\rm{(i) (Representation Theory)}}\enspace
$H \subset G \supset G'$ is a bounded multiplicity triple;
 that is, 
\[
\underset{\Pi \in {\text{\scriptsize{$\operatorname{Irr}(G)_H$}}}}\sup\,\,
\underset{\pi \in \operatorname{Irr}(G')}\sup
 [\Pi|_{G'}:\pi]< \infty.  
\]
{\rm{(ii) (Complex Geometry)}}\enspace
The generalized flag variety 
 $G_{\mathbb{C}}/B_{G/H}$ is $G_{\mathbb{C}}'$-spherical.  
\end{theorem}

\begin{example}
[$\operatorname{diag}G \nearrow G \times G \searrow G' \times G'$]
\label{ex:gGGg}
In view of the natural bijection, 
\begin{equation}
\label{eqn:GG}
   \operatorname{Irr}G 
   \simeq
   \operatorname{Irr}(G \times G)_{\operatorname{diag}G}, 
   \quad
   \pi \leftrightarrow \pi \boxtimes \pi^{\vee}, 
\end{equation}
 where $\pi^{\vee}$ is the contragredient representation, 
 Theorem \ref{thm:thm7} in this special case
 implies the bounded multiplicity theorem 
 \cite[Thm.\ D]{xktoshima}, 
 as stated in the equivalence (i) $\Leftrightarrow$ (ii) in Theorem \ref{thm:thm3}.  
\end{example}

The following theorem, 
 when  applied to the triple
 $\operatorname{diag}G \subset G \times G \supset G' \times G'$, 
 recovers the implication (ii) $\Rightarrow$ (i)
 in Theorem \ref{thm:thm4}
 ({\it{cf.}} \cite[Thm.\ C]{xktoshima}).

\begin{theorem}
[Finite Multiplicity Triple $H\nearrow G \searrow G'$ {\cite{K22}}]
\label{thm:thm8}
Retain the setting for $G' \subset G \supset H$.  
Then the following implications 
{\rm{(i)}} $\Rightarrow$ {\rm{(ii)}} $\Rightarrow$ {\rm{(iii)}} hold:
\newline
{\rm{(i)}} {\rm{(Complex geometry)}}\enspace
$\sharp (P_{\mathbb{C}}' \backslash G_{\mathbb{C}}/(P_{G/H})_{\mathbb{C}})<\infty$, 
 where $P'$ is a minimal parabolic subgroup of $G'$.  
\newline
{\rm{(ii) \enspace(Representation Theory)}}
The multiplicity $[\Pi|_{G'}:\pi] < \infty$ 
 for every $\Pi \in \operatorname{Irr}(G)_H$
 and for every $\pi \in \operatorname{Irr}(G')$.    
\\
{\rm{(iii) (Real Geometry)}}
The generalized real flag variety $G/{P_{G/H}}$ is $G'$-real spherical, 
 that is, 
 $P'$ has an open orbit in $G/{P_{G/H}}$
 (Definitions \ref{def:Rspherical} and \ref{def:PGH}).  
\end{theorem}

For the proof of Theorems \ref{thm:thm7} and \ref{thm:thm8}, 
 we prove a \lq\lq{QP estimate}\rq\rq\ 
 \cite[Thm.\ 3.1]{K22}
 for the uniform bounded multiplicity 
 of the restriction along the same line in \cite{xktoshima, Tu}
 and use a uniform upper estimate
 of the \lq\lq{largeness}\rq\rq\ 
 for all $H$-distinguished representations.  
The latter can be formulated via the following generalization
 of Harish-Chandra's subquotient theorem as outlined below.

\subsection{Generalizing Harish-Chandra's subquotient theorem}
~~~\newline

Harish-Chandra's subquotient theorem \cite{HC54b} has been strengthened
 as the {\it{subrepresentation theorem}}
 by Casselman \cite{CM82} among others:
 it asserts 
 that any $\pi \in \operatorname{Irr}(G)$
 can be realized
 as a subrepresentation
 (also as a quotient) of some principal series representation.

We present a theorem 
 that sharpen the subrepresentation theorem 
 for $\pi \in \operatorname{Irr}(G)_H$, 
 by replacing principal series representations
 with induced representations from finite-dimensional representations
 of the parabolic subgroup $P_{G/H}$.

In the special case of 
 $(G \times G, \operatorname{diag} G)$, 
 Theorem \ref{thm:GHquotient} recovers
 the subrepresentation theorem
 through the isomorphism \eqref{eqn:GG}.  
\begin{theorem}
[Subrepresentation Theorem for $G/H$]
\label{thm:GHquotient}
Let $(G,H)$ be a reductive symmetric pair.  
For any $\Pi \in \operatorname{Irr}(G)_H$, 
 there exists a finite-dimensional irreducible representation $V$ of {$P_{G/H}$} such that
\[
   \operatorname{Hom}_G(\operatorname{Ind}_{{\text{\scriptsize{$P_{G/H}$}}}}^G(V), \Pi)\ne \{0\}.  
\]
\end{theorem}

Moreover, 
 the representation $V$ in Theorem \ref{thm:GHquotient} has
 specific constraints
 that can be formulated
 using the relative Borel subalgebra ${\mathfrak{b}}_{G/H}$:
 see \cite[Thm.\ 1.8]{K22} for details.  
The proof of the uniform boundedness property in Theorem \ref{thm:thm7} 
 makes use of these constraints.

\subsection
{Classification: Bounded multiplicity triples}
~~~\newline
A complete classification of bounded multiplicity triples $H \subset G \supset G'$
 has been achieved in \cite{K21, K22};
 this classification is based on the criterion 
 in Theorem \ref{thm:thm7}.  
It has an unexpected relationship
 with the new spectral theory
 of locally symmetric spaces
 beyond the classical Riemannian setting, 
 which will be discussed
 in the next section.

To conclude this section, 
 we provide examples 
 of such bounded multiplicity triples.

\begin{example}
[$H \nearrow G \searrow G'$]
For any $p_1$, $q_1$, $p_2$, $q_2$
 with $p_1+p_2=p$, $q_1+q_2=q$, 
 the triple 
\[
(G, G', H)=(O(p,q), O(p_1,q_1) \times O(p_2,q_2), O(p-1, q))
\]
is a bounded multiplicity triple, 
 that is, 
\[
\underset{\text{\scriptsize{$\Pi \in \operatorname{Irr}(G)_H$}}}\sup
\,\,
\underset{{\text{\scriptsize{$\pi \in \operatorname{Irr}(G')$}}}}\sup
 [\Pi|_{G'}:\pi] < \infty.
\]
\end{example}

\section{Application of Branching Problem $G \downarrow G'$ in Geometry}
\label{sec:spec}
In this section, 
 we discuss a seemingly unrelated topic, 
 specifically the spectral analysis
 of pseudo-Riemannian locally symmetric spaces $\Gamma \backslash G/H$, 
 beyond the classical Riemannian setting.

We shall see that the general theory of the restriction $G \downarrow G'$, as discussed in Sections~\ref{sec:admrest} to \ref{sec:small}, provides a new tool for the study of $L^2(\Gamma \backslash G/H)$, where $H$ is a {\it non-compact} subgroup.

We begin with the general setup.  
Let $G$ be a Lie group, 
 $H$ be a closed subgroup, 
 and 
 $\Gamma$ be a {\it{discontinuous group}} for $X=G/H$.  
This means 
 that $\Gamma$ is a discrete subgroup of $G$
 acting properly discontinuously and freely on $G/H$.  
Consequently, 
 the double coset space $\Gamma \backslash G/H$, with the quotient topology, is Hausdorff
 and admits the unique $C^{\infty}$ manifold structure for which the covering map 
$p_{\Gamma} \colon G/H \to \Gamma \backslash G/H$ is a local diffeomorphism.  
$$
\begin{array}{@{}c@{}c@{}c@{}c@{}c@{}}
& & G  & & \\
& \swarrow & & \searrow & \\
\Gamma\backslash G  & & & & G /H  \\
& \searrow & & \swarrow & p_{\Gamma} \\
& & \Gamma \backslash G /H 
\end{array}
$$

\par\noindent
Via the covering map $p_{\Gamma}$, 
 any $G$-invariant local geometric object can be pushed
 forward to the quotient manifold $X_{\Gamma}:=\Gamma \backslash G/ H$.

Let ${\mathbb{D}}_G(X)$ denote
 the algebra of $G$-invariant differential operators on $X=G/H$.  
Then any $D \in {\mathbb{D}}_G(X)$ induces 
 a differential operator $D_{\Gamma}$
 on the quotient $X_{\Gamma}$
 via the covering map $X \to X_{\Gamma}$.  
We consider the set
\[
   {\mathbb{D}}(X_{\Gamma}):=\{D_{\Gamma}: D \in  {\mathbb{D}}(G/H)\}
\]
 as the algebra
 of {\it{intrinsic differential operators}}
 on the locally homogeneous space $X_{\Gamma}$.

\begin{example}
If $G \supset H$ is a pair of real reductive Lie groups, 
 then $X_{\Gamma}$ inherits a pseudo-Riemannian structure from a $G$-invariant pseudo-Riemannian structure on $X=G/H$.  
The Laplacian $\Delta_{X_{\Gamma}}$
 is an element of ${\mathbb{D}}(X_{\Gamma})$.  
\end{example}

We remark that for non-compact $H$, 
 the pseudo-Riemannian structure is not necessarily positive definite; consequently, 
 the Laplacian is generally not an elliptic differential operator.  
For instance, 
 if the quotient space $G/H$ is Lorentzian, 
 the following equation
 on the space $X_{\Gamma}$:
\begin{equation}
\label{eqn:Laplmd}
  \Delta_{X_{\Gamma}} f = \lambda f 
\end{equation}
is a hyperbolic equation.

Suppose now that $X=G/H$ is a reductive symmetric space.  
Then the algebra 
$
   {\mathbb{D}}(X_{\Gamma})\simeq {\mathbb{D}}_G(X)
$
 is a commutative ring.  
Hence, 
 it is natural to consider joint eigenfunctions
 of ${\mathbb{D}}(X_{\Gamma})$
 on $C^{\infty}(X_{\Gamma})$
 rather than focusing on the single equation \eqref{eqn:Laplmd}.

For any ${\mathbb{C}}$-algebra homomorphism
$
   \lambda \colon {\mathbb{D}}(X_{\Gamma}) \to {\mathbb{C}}, 
$
we denote by $C^{\infty}(X_{\Gamma}, {\mathcal{M}}_{\lambda})$
 the space of smooth functions
 defined on $X_{\Gamma}$ that satisfy the system of equations
\begin{equation}
  D_{\Gamma} f = \lambda(D) f \quad
  \text{for any } D \in {\mathbb{D}}_G(X).  
\tag{${\mathcal{M}}_{\lambda}$}
\end{equation}
Let ${\mathbb{D}}_G(X)^\wedge$ be the set of 
 all ${\mathbb{C}}$-algebra homomorphisms
 $\lambda \colon {\mathbb{D}}(X_{\Gamma}) \to {\mathbb{C}}$, 
 which can be identified with the quotient
 ${\mathfrak{j}}_{\mathbb{C}}^{\vee}/W$
 of the dual ${\mathfrak{j}}_{\mathbb{C}}^{\vee}$
 of a Cartan subspace ${\mathfrak{j}}$
 for $(G, H)$
 by the Weyl group $W$
 for the root system 
$
   \Delta({\mathfrak{g}}_{\mathbb{C}}, {\mathfrak{j}}_{\mathbb{C}})
$, 
 as shown in the Harish-Chandra isomorphism
 \cite{HC8}.

Not much attention has been paid to the spectral theory 
 on $X_{\Gamma}=\Gamma \backslash G/H$ 
 in the general setting where $H$ is {\it{non-compact}}
 and $\Gamma$ is an infinite discontinuous group.  
For instance, 
 the following questions regarding the spectral theory
 remain open in this generality.  

\begin{prob}
\label{prob:specGamma}
\item[{\rm{(1)}}]
The expansion of arbitrary functions 
 defined on $X_{\Gamma}$
 in terms of joint eigenfunctions
 of the algebra ${\mathbb{D}}(X_{\Gamma})$
 of intrinsic differential operators.  
\item[{\rm{(2)}}]
Understanding the distributions of $L^2$-eigenvalues.  
\end{prob}


\subsection{Spectral analysis on $\Gamma \backslash G /H $ in the classical case}
~~~\newline
These problems for spectral analysis 
 on $X_{\Gamma}=\Gamma\backslash G/H$
 are formulated from a broader perspective, 
 building on the rich and deep body 
 of classical results
 that have been extensively studied.  
Special (classical) cases 
 that have been particularly fruitful include:

\begin{enumerate}
\item[(1)]
Let $H=K$.  
When $H$ is a maximal compact subgroup $K$ of $G$, 
 $X_{\Gamma}$ becomes a {\it{Riemannian}} locally symmetric space.  
A vast theory has been developed
 over several decades.  
Problem \ref{prob:specGamma} is particularly enriched
 in connection with the local theory 
 of automorphic forms
 when $\Gamma$ is an arithmetic subgroup.  

\item[(2)]
Let $\Gamma=\{e\}$.  
When $\Gamma=\{e\}$, $X_{\Gamma}$ reduces to the homogeneous space $G/H$.  
Problem \ref{prob:specGamma} (1) has been extensively studied in the case
 where $G/H$ is a reductive symmetric space
 and $\Gamma=\{e\}$, 
 with significant contributions from Gelfand,
 Harish-Chandra \cite{HC76}
 for group manifolds $(G \times G)/\operatorname{diag}G$;
 Helgason, Flensted-Jensen, T.\,Oshima, Delorme, 
 and others for symmetric spaces $G/H$.  
\item[(3)]
Let $G={\mathbb{R}}^{p,q}$, 
 $\Gamma={\mathbb{Z}}^{p+q}$, 
 and $H=\{0\}$.  
In this case, 
 $\operatorname{Spec}_d(X_{\Gamma})$ consists 
 of values of indefinite quadratic forms
 of signature $(p,q)$
 at the dual lattice $\Gamma^{\vee}$.  
See \cite{K16a}
 for a discussion on Problem \ref{prob:specGamma} (2), 
 which highlights a relation 
 to the Oppenheim conjecture
 (see {\it{e.g.}}, Margulis \cite{Mg00} and references therein)
 in Diophantine approximation.  
\end{enumerate}

\subsection{Spectral analysis on $\Gamma \backslash G/H$
 beyond the classical Riemannian setting}
\label{subsec:challenge}
~~~\newline
The situation changes drastically
 beyond the classical setting, 
 that is, 
 when $H$ is no longer compact and $\Gamma \ne \{e\}$.  
New difficulties arise, 
including:

\begin{enumerate}
\item[(1)]
(Representation theory)\enspace
The space $L^2(X_{\Gamma})$ is no longer a subspace of $L^2(\Gamma \backslash G)$, 
 where the group $G$ acts as a unitary representation.  
Moreover, 
 we cannot expect the regular representation of $G$ 
 on $L^2(\Gamma \backslash G)$
 to have finite multiplicities, 
 contrary to the classical theorem 
 of Gelfand--Piatetski--Shapiro.  
\item[(2)]
(Analysis)\enspace
In contrast to the Riemannian case, 
 where $H=K$, 
 the pseudo-Riemannian Laplacian $\Delta_{X_{\Gamma}}$ is
 no longer an elliptic differential operator.  
Moreover, 
 it is unclear 
 whether $\Delta_{X_{\Gamma}}$ is essentially self-adjoint, 
 due to the absence of a general theory.  

\item[(3)]
(Geometry)\enspace
If $H$ is not compact, 
 then not all homogeneous spaces $G/H$ can admit discontinuous groups
 of infinite order
 ({\it{e.g.,}} the Calabi Markus phenomenon \cite{CM62, kob89}).  
\end{enumerate}


\subsection{Standard locally homogeneous spaces $\Gamma\backslash G/H$}
~~~\newline
The geometric issue (3) in Section \ref{subsec:challenge}
 raises a \lq\lq{local to global}\rq\rq\ problem
 beyond the Riemannian setting, 
 which can be formulated
 by group-theoretic terms as the following fundamental question.

\begin{question}
\label{q:discont}
How can we find a discrete subgroup $\Gamma$
 that acts properly discontinuously on $G/H$?
\end{question}

Here are some elementary observations 
regarding two extreme cases.

\begin{observation}
(1)\enspace
Any discrete subgroup $\Gamma$ suffices 
 if $H$ is compact.  
\newline
(2)\enspace
However, 
 any lattice $\Gamma$ of the whole group $G$ does not suffice 
 if $H$ is non-compact, 
 because $\Gamma \backslash G/H$ cannot be Hausdorff
 in this case.  
\end{observation}

A straightforward method
 to find discrete subgroups $\Gamma$
 that answer Question \ref{q:discont} is
 to utilize \lq\lq{standard quotients}\rq\rq:
we recall
 that $G'$ acts {\it{properly}}
 if the map $G' \times X \to X \times X$, 
 given by $(g, x) \mapsto (x,gx)$, 
 is proper.

\begin{definition}
\label{def:standard}
Let $G'$ be a reductive subgroup, 
 which acts properly on $X$, 
 and let $\Gamma$ be a torsion-free discrete subgroup of $G'$.  
Then $\Gamma$ acts properly discontinuously on $X$, 
 and we say the quotient manifold $X_{\Gamma}= \Gamma \backslash G/H$ is 
 a {\it{standard quotient}} of $X$.  
\end{definition}

By using standard quotients, 
 it turns out
 that there exist several families of reductive symmetric spaces $G/H$
 that admit \lq\lq{large}\rq\rq\ discontinuous groups $\Gamma$, 
 {\it{e.g.,}} such that $X_{\Gamma}=\Gamma \backslash G/H$ is compact 
or has finite volume
 \cite{Ku81,kob89}.


\subsection{Spectral theory of standard locally symmetric space 
$\Gamma\backslash G/H$}
~~~\newline
For the study of spectral analysis
 of pseudo-Riemannian locally symmetric spaces, 
 we focus on {\it{standard}} quotients $X_{\Gamma}$
 (Definition \ref{def:standard})
 of a reductive symmetric space $X=G/H$, 
 where $\Gamma$ is a discrete subgroup
 of a reductive subgroup $G'$ 
 that acts properly on $X$.  

\begin{theorem}
[Expansion into Joint Eigenfunctions, {\cite{KK20, SpecLNM2024}}]
\label{thm:thm9}
Assume that $G_{\mathbb{C}}'$ acts 
spherically on $X_{\mathbb{C}}$
 (Definition \ref{def:spherical}).  
Then any $f \in C_c^{\infty}(\Gamma\backslash X)$ can be expanded into 
 joint eigenfunctions of ${\mathbb{D}}(X_{\Gamma})$ on $\Gamma\backslash X$.  
More precisely, 
 there exist a measure $\mu$
 on ${\mathbb{D}}_G(X)^\wedge (\simeq {\mathfrak{j}}_{\mathbb{C}}^{\vee}/W)$
 and
 a family of maps
\\
\[
  {\mathcal{F}}_\lambda \colon 
  C_c^{\infty}(X_{\Gamma}) \to C^{\infty}(X_{\Gamma};{\mathcal{M}}_{\lambda})
\]
such that 
\[
f = \displaystyle\int_{{\mathbb{D}}_G(X)^\wedge} {\mathcal{F}}_{\lambda} f \, d \mu(\lambda), 
\]
\[
\|f\|_{L^2(X_{\Gamma})}^{2} = \int_{{\mathbb{D}}_G(X)^\wedge}
\|{\mathcal{F}}_{\lambda}f\|_{L^2(X_{\Gamma})}^2 d \mu(\lambda)
\]
 for any $f \in C_c^{\infty}(X_{\Gamma})$.  
\end{theorem}

\begin{remark}
In Theorem \ref{thm:thm9}, 
 we do not assume
 that $\operatorname{vol}(\Gamma \backslash G') < \infty$.  
In particular, 
 Theorem \ref{thm:thm9} holds even when $\Gamma=\{e\}$.  
\end{remark}

The following corollary answers an analytic issue
 (2) in Section \ref{subsec:challenge}
 in the setting of Theorem \ref{thm:thm9}.  

\begin{corollary}
\label{cor:Lap}
In the setting of Theorem \ref{thm:thm9}, 
 the pseudo-Riemannian Laplacian $\Delta_{X_{\Gamma}}$
 is essentially self-adjoint on $L^2(X_{\Gamma})$.  
\end{corollary}


\subsection{Examples}
~~~\newline
The properness criterion for the triple $(G', G, H)$
 was established in \cite{kob89}, 
 and is both explicit and computable.  
Below are some examples
 to which Theorem \ref{thm:thm9} is applied.

\begin{example}
[Riemannian locally symmetric space]
Let $G/K$ be a Riemannian symmetric space.  
Set $G':=G$.  
Then the triple $(G', G, K)$ clearly satisfies
 the assumption of Theorem \ref{thm:thm9}.  
Consequently, 
 the conclusion of Theorem \ref{thm:thm9} holds
 for any Riemannian locally symmetric space 
 $\Gamma \backslash G/K$, 
 including the case of infinite volume.  
\end{example}

\begin{example}
[standard anti-de Sitter manifolds]
Let $X$ be an odd-dimensional 
 anti-de Sitter space, 
 that is, 
 $X=G/H=SO(2n,2)/SO(2n,1)$.   
The subgroup $G':=U(n,1)$ acts properly on $X$, 
 and $G_{\mathbb{C}}'=G L(n+1, {\mathbb{C}})$ acts spherically
 on $X_{\mathbb{C}}=SO(2n+2, {\mathbb{C}})/SO(2n+1, {\mathbb{C}})$.  
Therefore, 
 Theorem \ref{thm:thm9} holds for any standard anti-de Sitter manifold
 $X_{\Gamma}$ 
 with $\Gamma \subset U(n,1)$.  
\end{example}

\begin{example}
[indefinite K{\"a}hler manifolds]
The homogeneous space $X=G/H=SO(2n,2)/U(n,1)$ has
 a natural indefinite-K{\"a}hler structure, 
 and the subgroup $G'=SO(2n,1)$ acts properly.  
Moreover, 
 $G_{\mathbb{C}}'=S O(2n+1, {\mathbb{C}})$ acts spherically
 on $X_{\mathbb{C}}= S O(2n+2, {\mathbb{C}})/G L(n+1,{\mathbb{C}})$.  
Thus, 
 Theorem \ref{thm:thm9} holds
 for any standard indefinite-K{\"a}hler manifold $X_{\Gamma}$
 with $\Gamma \subset S O(2n,1)$.  
\end{example}

\begin{example}
[15-dimensional space form of signature $(8,7)$]
Let $X=G/H=SO(8,8)/SO(8,7)$.  
Then $X$ is a pseudo-Riemannian manifold of signature $(8, 7)$
 with constant negative sectional curvature.  
The subgroup $G'=Spin(7,1)$ of $G$ acts properly on $X$, 
 and $G_{\mathbb{C}}'= S p i n(8, {\mathbb{C}})$
 acts spherically on $X_{\mathbb{C}}$.  
Thus, 
 Theorem \ref{thm:thm9} holds
 for any standard 15-dimensional space from $X_{\Gamma}$
 of signature (8,7)
 with $\Gamma \subset S p i n (7,1)$.  
\end{example}

A classification of the triples
 $G' \subset G \supset H$ 
 that satisfy
 the assumptions in Theorem~\ref{thm:thm10}---specifically, 
 the subgroup $G'$ acts properly on $X=G/H$
 and $X_{\mathbb{C}}$ is $G_{\mathbb{C}}'$-spherical-- 
 can be found in \cite{KK20, SpecLNM2024}.   
As will be observed in Section \ref{subsec:5.7}, 
 such triples are {\it{bounded multiplicity triples}}
 (Definition \ref{def:3.3}) discussed in the previous section.


\subsection{Key step: Admissible restriction $H \nearrow G \searrow G' \supset \Gamma$}
~~~\newline
The new approach for proving the spectral theory
 presented in Theorem \ref{thm:thm9}
 involves utilizing the global analysis of $G_{\mathbb{C}}'$-spherical spaces
 and the restriction of irreducible representations $\Pi$ of $G$
 to the subgroup $G'$
 that acts properly on $X$.  
Specifically, 
 we demonstrate
 that this restriction $\Pi|_{G'}$ is always
 discretely decomposable
 with uniformly bounded multiplicities.

The following theorem bridges the spectral analysis
 on $X_{\Gamma}=\Gamma \backslash G/H$
 and recent progress of branching problems
 that we outlined in Sections \ref{sec:admrest} to \ref{sec:small}.  

\begin{theorem}
\label{thm:thm10}
Let $X=G/H$ be a reductive symmetric space.  
Suppose that a reductive subgroup $G'$ of $G$ 
 acts properly on $X$ 
and that the complexification $G_{\mathbb{C}}'$ acts spherically
 on $X_{\mathbb{C}}$.  
Then the restriction $\Pi|_{{\text{\tiny{$G'$}}}}$
 is {$G'$-admissible}
 (Definition \ref{def:deco})
 for any $H$-distinguished $\Pi \in \widehat G$
 (Definition \ref{def:distinguished}).  
Moreover, 
 the multiplicities are uniformly bounded:
\[
  \underset{\Pi \in \operatorname{I r r}(G)_H}\sup
\,\,
  \underset{\pi \in \operatorname{I r r}(G')}\sup
  [\Pi|_{G'}:\pi]<\infty.  
\]
\end{theorem}


\subsection{
$H \,\,\subset\,\, G \,\,\supset\,\, G'
\,\,\,\,\rightsquigarrow\,\,\,\,\,\,
H\,\,
{}^{\text{\tiny{ind}}}\hskip -0.8pc\nearrow 
G
\searrow
\hskip -0.6pc
{}^{\text{\tiny{rest}}}
\,\,
G'
$
}
\label{subsec:5.7}
~~~\newline
The conclusion of Theorem \ref{thm:thm10} follows from 
 two key aspects
 of branching problems:
 $G'$-admissibility, 
 as discussed in Section \ref{sec:admrest}
 and the uniformly bounded multiplicities
 outlined in Sections \ref{sec:realsp} and \ref{sec:small}.

The following diagram summarizes the related results:
\begin{alignat*}{3}
&G_{\mathbb{C}}' 
\underset{\text{spherical}}{\curvearrowright} 
G_{\mathbb{C}}/B_{G_{\mathbb{C}}/H_{\mathbb{C}}}
&&\overset{\text{Theorem \ref{thm:thm7}}}\iff\quad
&&\underset{\Pi \in \operatorname{Irr}(G)_H}{\sup}\,\,
\underset{\pi \in \operatorname{Irr}(G')}{\sup}
[\Pi|_{G'}:\pi]<\infty
\\[1em]
&\hphantom{MMMM} \rotatebox[origin]{90}{$\Longrightarrow$}
&& 
&&\hphantom{MMMMMM} \rotatebox[origin]{90}{$\Longrightarrow$}
\\[1em]
&
\begin{cases}
G_{\mathbb{C}}' 
\underset{\text{spherical}}{\curvearrowright} 
G_{\mathbb{C}}/H_{\mathbb{C}}
\\
G' 
\hphantom{i}\underset{\text{proper}}{\curvearrowright}
G/H
\end{cases}
&&\underset{\text{Theorem \ref{thm:thm10}}}{\Longrightarrow}\quad
&& 
{\begin{aligned}
&\text{$\Pi|_{G'}$ is {discretely decomposable}}
\\[-.3em]
&\text{with {uniformly bounded multiplicities}}
\\[-.3em]
&\text{for all $\Pi \in \operatorname{Irr}(G)_H$.}
\end{aligned}
}
\end{alignat*}

The proof of Theorem \ref{thm:thm10}
 (\cite{K16b}, {\it{cf.}} \cite{KK20, SpecLNM2024})
  does not rely
 on the $G'$-admissibility criterion
 (Theorem \ref{thm:thm1}), 
 which is formulated purely
 in terms of representation theory. 
It would be interesting 
 to explore a direct connection 
 between the geometric assumption 
 in Theorem \ref{thm:thm10}
 and the transversality condition 
 of the two cones
 (condition (ii) in Theorem \ref{thm:thm1}).


\subsection{Strategy of the proof for Theorem \ref{thm:thm9}}
~~~\newline
The proof of Theorem \ref{thm:thm9}
 in \cite{KK20, SpecLNM2024}
 is quite comprehensive.  
Here, 
 we outline the key ingredients.  
Recall that we consider the standard quotient:
\[
  \Gamma \subset 
G' \subset G 
{\raisedcurvearrowright}
 X
\quad
  \rightsquigarrow 
\quad
X_{\Gamma}=\Gamma \backslash X
\]
\begin{enumerate}
\item[1.]
(Hidden symmetry)
If the action of $G_{\mathbb{C}}'$ on $X_{\mathbb{C}}$ is spherical 
 (see Theorem \ref{thm:thm2}), 
 it can be shown 
 that the algebra ${\mathbb{D}}_{G'}(X)$ leaves
 the space $C^{\infty}(X;{\mathcal{M}}_{\lambda})$ invariant, 
 see \cite{K16b, KaTK19}.  
In other words, 
 ${\mathbb{D}}_{G'}(X)$ acts as a hidden symmetry 
 of the joint eigenspace $C^{\infty}(X; {\mathcal{M}}_{\lambda})$:
\[
  {\mathbb{D}}_{G}(X)
 \subset {\mathbb{D}}_{G'}(X) \raisedcurvearrowright C^{\infty}(X;{\mathcal{M}}_{\lambda}).  
\]
\item[2.]
(Branching law $G \downarrow G'$)
If the $G'$-action on $X$ is also proper, 
 then any $\pi \in \operatorname{Irr}(G)$ realized
 in ${\mathcal{D}}'(X)$ is 
 {$G'$}-admissible
 (Theorem \ref{thm:thm10}).  
This property is particularly clear
 in the special case where $H=K$ and $G'=G$, 
 as it 
 corresponds to Harish-Chandra's admissibility.  
\end{enumerate}

\end{document}